\documentclass[10pt, a4paper, reqno]{amsart}

\usepackage{amsmath, amssymb, amsfonts, enumerate}
\usepackage[english]{babel}
\usepackage{fancyhdr}
\usepackage{graphicx}
\usepackage{amsfonts}
\usepackage{mathrsfs}
\usepackage{mathptmx}
\usepackage{stmaryrd}
\usepackage{hyperref}

\theoremstyle{plain}

\newtheorem*{teo*}{Theorem}
\newtheorem{teo}{Theorem}
\newtheorem{rem}{Remark}

\newcommand{\N}{\mathbb{N}}
\newcommand{\F}{\mathbb{F}}
\newcommand{\Z}{\mathbb{Z}}
\newcommand{\T}{\mathbb{T}}

\newcommand{\R}{\mathbb{R}}

\makeatletter
\def\imod#1{\allowbreak\mkern10mu({\operator@font mod}\,\,#1)}
\makeatother

\newcommand{\dig}[2]{\llbracket #1 \rrbracket_{#2}}

\begin{document}

\title[On multipliers providing perfect hashing]{On perfect hashing
 of numbers with sparse digit representation
via multiplication by a constant}

\author[M. Monge]{Maurizio Monge}
\address{Scuola Normale Superiore di Pisa - Piazza dei Cavalieri, 7 - 56126 Pisa}
\email{maurizio.monge@sns.it}
\date{\today}

\begin{abstract}
  Consider the set of vectors over a field having non-zero
  coefficients only in a fixed sparse set and multiplication defined
  by convolution, or the set of integers having non-zero digits (in
  some base $b$) in a fixed sparse set. We show the existence of an
  optimal (resp. almost-optimal in the latter case) `magic' multiplier
  constant that provides a perfect hash function which transfers the
  information from the given sparse coefficients into consecutive
  digits. Studying the convolution case we also obtain a result of
  non-degeneracy for Schur functions as polynomials in the elementary
  symmetric functions in positive characteristic.
\end{abstract}

\maketitle

\section{Introduction and Motivation}
Suppose $n>0$, and let $D = \{ d_0, d_1, \dots, d_{n-1} \}$ be a set
of indices such that $0 = d_0 < d_1 < \dots < d_{n-1}$. For a field
$F$ and $N = d_{n-1}+1$, let $F^{N}$ be the standard vector space
with basis $e_j$ for $0 \leq j < N$ equipped with
convolution multiplication
\[ (a \ast b)_i = \sum_{\substack{j+k=i\\0\leq j,k < N}} a_j b_k, \]
 and let $F[D]$ be the subspace spanned by
$e_{d_0},\dots,e_{d_{n-1}}$, which is formed by the vectors with
non-zero coefficients only in the indices $d_i$ for $0\leq i \leq
n-1$.  Similarly given a positive integer $b$ let $\Z_b[D]$
be the set of integers that can be written as $\sum_{i=0}^{n-1} a_i
b^{d_i}$ for some $a_i \in \{0,1,\dots, b-1\}$, i.e. the set of
numbers such that their base-$b$ representation only contains
non-zero digits in positions that belong to the set $D$.  We study the
existence of constants that can be used as multipliers to transfer the
information stored in the sparse digits of an element of $F[D]$ or of
$\Z_b[D]$ into a smallest possible set of consecutive digits,
providing a perfect hash function.

Motivation for this kind of questions is provided by a technique used
by many state-of-the-art chess playing programs \cite{pradu2007magic,
  tannous2007avoiding}, based on the concept of `bitboards', that are
numbers whose base-$2$ representation is interpreted as an occupancy
information of some kind, or more generally to store a $0-1$
information for each square, having previously established a
correspondence between a range of digits and the squares on the
board. The technique in question, known under the name of `magic
bitboards', is a quick way to generate all possible attacks for
sliding pieces, such as rooks and bishops. The bitboard containing
occupancy information for all pieces is transformed with a bitwise-and
to a bitboard whose only digits that may be different from zero are
those corresponding to possible obstructions on the path of the
sliding piece. This information about the obstructions, which is
stored in a small set of sparse digits, is then mapped via a
multiplication by a `magic number' to a set of consecutive digits,
which is then used as index in a lookup table to recover a
pre-calculated information about the possible attacks.

While in the case of chess programs a database of very efficient
multipliers has already been computed and is publicly available, we
investigate the existence of multipliers that provide perfect hashing
functions in a more general setting.  In the convolution case we
provide an optimal result, which shows that it is possible to transfer
the information stored in any number of sparse digits into the same
number of consecutive digits, and which incidentally provides a result
about values of Schur functions as polynomials in elementary symmetric
functions.  On the other hand, in the case of base-$b$ integers it is
not always possible to have a multiplier providing a hash into the
same number of digits ($D=\{0,1,2,4,6\}$ providing the smallest
counterexample for $b=2$, as can be checked with a simple
computer program), and we provide a linear estimate of the number of
consecutive digits that are required to ensure the existence of such a
map.

While this kind of hashing cannot be directly compared to universal
hashing (see \cite{thorup2000even,thorup2004tabulation,
  raman1996priority}, and \cite{woelfel1999efficient} in particular)
because of its more restricted scope, it is still possible to compare the
results about its effectiveness, and this is done below. See also
\cite{sauerhoff2003time}.

\subsection*{Acknowledgements} We wish to thank Vincenzo Mantova for
the time enjoyed discussing this and related questions. We also thank
the reviewers for suggesting relevant references.

\section{The convolution case}

In this section we consider the convolution case. The operation of
taking the convolution multiplication with a fixed vector
$(a_0,\dots,a_{N-1})$ can be expressed by a lower triangular Toeplitz
matrix $A = (a_{i-j})_{0\leq i,j < N}$, where we have put $a_i = 0$
for negative $i$ (the full convolution with a vector
$(a_{-N+1},\dots,a_{N-1})$ is expressed by a general Toeplitz matrix, but as
shown below we can restrict to the class of lower triangular
matrices). In the same way, when restricting the output to a set of
coefficients we shall consider the matrix formed by the corresponding
selection of rows of $A$.  We will show that for a good choice of the
$a_i$ the matrix formed by the last $n$ rows of $A$ defines a
one-to-one function from $F[D]$ to $F^n$.

The operation of taking the convolution with a fixed vector, or multiplying by a Toeplitz matrix, is
also equivalent to the multiplication by a polynomial in a ring of polynomials, and its properties
as hash function are well known, as well as fast algorithms, see
\cite{mansour1990computational,gohberg1994complexity}. Note also that taking the convolution with a
vector with entries in the set $\{0,1\}$ is actually obtained by addition of selected entries.

\begin{teo} \label{teo1}
  For every set $D$ of cardinality $n$, there exists a lower triangular Toeplitz matrix such that
  its last $n$ rows define a one-to-one function from $F[D]$ to $F^n$, and furthermore its entries
  can be taken in the set $\{0,1\}$.
\end{teo}

\begin{proof}
  Let $\delta = N-n$, and for a vector $(a_0,\dots,a_{N-1})$ that is going to be detemined let $A$
  be the associated rectangular Toeplitz matrix $(a_{\delta+i-j})_{\substack{0\leq i < n\\0\leq j < N}}$
  corresponding to the selection of the last $n$ coefficients after taking the convolution. Let $B$
  be the $n\times n$ matrix obtained from $A$ selecting the columns with indices
  $d_0,\dots,d_{n-1}$, which defines the linear map on the basis of $F[D]$ formed by the $e_{d_i}$
\[
   B = (a_{\delta+i-d_k})_{\substack{0\leq i
      < n\\0\leq k < n}}
 =
  \begin{pmatrix}
 a_{\delta-d_0} & a_{\delta-d_1} & \dots & a_{\delta-d_{n-1}} \\ 
 a_{\delta+1-d_0} & a_{\delta+1-d_1} & \dots & a_{\delta+1-d_{n-1}} \\
 \vdots  & \vdots  &          & \vdots    \\
 a_{\delta+n-1-d_0} & a_{\delta+n-1-d_1} & \dots & a_{\delta+n-1-d_{n-1}} \\
  \end{pmatrix}.
  \]
  We will now consider a sequence of $k\times k$ minors, for
  $k=1,\dots,n$, where each minor will contain the previous one,
  and inductively change some of the $a_i$ ensuring at the $k$-th step that
  the determinant of the $k$-th minor is non-zero, while leaving
  unchanged the coefficients of the minors considered in the previous
  steps.

  Let the $k$-minor $B_k$ be obtained taking the first $k$ rows, and a
  range of columns $r_k,r_k+1,\dots,r_k+k-1$, where $0\leq r_k \leq
  n-k$ is the biggest integer such that each column of $B_k$ will
  contain at least one $a_i$ with $i\geq 0$. It is an easy consequence of $B$ being a
  selection of columns from the Toeplitz matrix $A$ that
  the set of integers $r_k$ satisfying the above condition is always
  non-empty, and that the set of columns selected for $B_k$ is the
  same set that was selected for $B_{k-1}$ with one column added
  either on the left or on the right.

  Let's begin the induction putting $a_0 = 1$ and $a_i = 0$ for all $i
  \neq 0$. For $k=1$, change also $a_{\delta-d_{r_1}}$ to be equal to
  $1$.  Let now $k$ be $>1$. When the columns of $B_k$ are those of
  $B_{k-1}$ plus one column on the right, we have that the matrix
  $B_k =
  \begin{pmatrix}
    B_{k-1} & 0 \\ \ast & a_0
  \end{pmatrix}$ is block lower triangular, with one block equal to
  $B_{k-1}$, and the other block being formed by the element $a_0 =
  1$, and $B_k$ is non-singular.  On the other hand, when one column
  is added on the left, $B_k$ is of the form $B_k =
  \begin{pmatrix}
    \ast & B_{k-1} \\ a_\ell & \ast
  \end{pmatrix}$, where $\ell = \delta+k+1-d_{r_k}$ is the biggest index appearing in $B_k$ (in
  fact, the indices are decreasing while moving right along a row or up along a
  column). Consequently, considering the Laplace expansion of the determinant of $B_k$ along the
  first column
\begin{equation} \label{eq1.1}
   \det B_k = \sum_{i=1}^{k-1} a_{\ell-k+i} \cdot C_{i,1}(k) + a_\ell \cdot \det B_{k-1}
\end{equation}
 where for each $i,j$, $C_{i,j}(k)$ is the $i,j$ cofactor of the matrix $B_k$,
 we can select an appropriate value for $a_\ell$ which makes the determinant non-zero, while changing only
  the bottom left entry of $B_k$. Repeating this step up to $k=n$ we have the theorem.
\end{proof}

\begin{rem}
  Since no division is involved the above proof works in any ring with $1$, but it does not ensure
  the resulting matrix $B$ to be invertible, only to have non-zero determinant. Alternatively, it is
  possible to allow general $a_i$, and solve inductively the \eqref{eq1.1} in $a_\ell$ to ensure
  that $\det B_k = 1$ at each step, obtaining that $B$ is invertible because the determinant is an
  invertible element of the ring.

  If $F$ is a local ring (i.e. a ring with only one maximal ideal) with maximal ideal $M$ we can
  still take the $a_i$ in $\{0,1\}$, and make all the determinants of the $B_k$ invertible: consider
  the \eqref{eq1.1} and call $S$ the sum, since in a local ring the invertible elements are
  precisely those not in $M$ and we assume $\det B_{k-1}$ to be invertible, we cannot have both
  $S\in{}M$ and $S+\det{}B_{k-1}{}\in M$ or we would also have $\det B_{k-1} \in M$, and hence
  putting $a_\ell$ equal to either $0$ or $1$ we obtain that $\det{}B_k\not\in{}M$, i.e. that it is
  invertible.

  This has a practical consequence: take for instance $F$ to be the set of integers modulo $2^k$, it
  is a local ring and hence we have the existence of a magic multiplier with entries in the set
  $\{0,1\}$, and taking the convolution is actually the addition of selected enties.
\end{rem}

It is possible to observe that the Schur functions that are well known
in Algebraic Combinatorics (see \cite{macdonald1995,
  sagan2001symmetric}) for their combinatorial properties and
connections with the characters of the symmetric group can be
expressed as a determinant of a special matrix having the elementary
symmetric functions as coefficients via the Jacobi-Trudi identity
(also known as ``determinant formula''), which has the same form as
the transpose of the matrix $B$ considered above. In particular, if
$s_\lambda$ is the Schur function associated to the partition $\lambda
= (\lambda_1,\lambda_2,\dots)$ and $\lambda'$ is the conjugate
partition, and $e_i$ is the $i$-th elementary symmetric function for
$i\geq 0$, we have the formula
\[ s_\lambda = det(e_{\lambda'_i-i+j})_{1\leq{}i,j\leq{}n} \] 
expressing the Schur function $s_\lambda$ as a uniquely determined polynomial in the $e_i$ having
integral coefficients. Since the set $D$ considered in Theorem \ref{teo1} is arbitrary, for each
partition of length $n$ we chose $D$ ensuring that $\delta-d_{k-1} = \lambda_k'-k$ for
$1\leq{}k\leq{}n$, making the matrix $B$ equal to the transposed of above matrix evaluated with
$e_i=a_i$ for all $i$. It follows that the polynomials expressing the $s_\lambda$ in terms of the
$e_i$ assume a non-zero value when the $e_i$ are replaced with opportune values in the field $F$
(note that this is not true for a general polynomial over a field, as the example $x^p-x$ over
$\F_p$ shows). Hence we have

\begin{teo}
  For a partition $\lambda$, consider the Schur function $s_\lambda$ as a polynomial in the
  elementary symmetric functions $e_i$, considered as indeterminates. Then it takes a non-zero value
  after substitution of the variables $e_i$ with appropriate elements of the field $F$, which
  moreover can be taken in the set $\{0,1\}$.
\end{teo}

\section{The arithmetic case}

The case of base-$b$ digits of integers seems to be much more
difficult, and we give a linear bound on the number of digits required
to ensure the existence of an opportune multiplier. Let $\Z_{(b)}$ be
the set of rational numbers that can be written as $r/b^k$ for some
integers $r,k$, or equivalently that have a finite base-$b$
expansion. For such an $a = \sum_{i<M} a_ib^i \in \Z_{(b)}$, define
$\dig{a}{k,m}$ as the $m$-tuple
$(a_{k+m-1},\dots,a_{k+1},a_k)\in\{0,1,\dots,b-1\}^m$. We can now state

\begin{teo}
Let $D = \{d_0,\dots,d_{n-1}\}$ be a set of indices $0 = d_0 < \dots < d_{n-1}$ having
cardinality $n$. Then for
\[
 m = \left\lceil \log_b\left((2b-1)^n-1\right) \right\rceil
\]
there exist a multiplier $\mu \in \Z$ and a $k \in \N$ such that the map from $\Z_b[D]$ to
$\{0,1,\dots,b-1\}^m$ defined by $a \mapsto \dig{a \cdot \mu}{k,m}$ is
injective. Furthermore, if $D$ contains some consecutive integers and is
formed by the union of the integral intervals $\{c_i,c_i+1,\dots,c_i+\ell_i-1\}$ for
$i=1,\dots,k$ and $\ell_i\geq 1$, we can take
\[
 m = \left\lceil \log_b\left( \prod_{i=1}^k (2b^{\ell_i}-1)-1 \right)\right\rceil.
\]
\end{teo}

It is possible to compare this estimate with what can be obtained
using universal hashing: when a hashing function is randomly chosen in
a universal class (i.e. $h(x) = h(y)$ with probability at most
$1/b^m$, what can be done when the output is formed by at least $m$
digits, see \cite{woelfel1999efficient}) we have that the probability
of $h$ being one-to-one on $\Z_b[D]$ is at least $1 -
\binom{b^n}{2}b^{-m}$, and we deduce the existence of a good hash
function when $m$ is $\geq 2n$. The above result is sharper because it
just requires $m$ to be about $n\log_b(2b-1)$, which is always smaller
than $2n$ and its ratio with $n$ approaches $1$ as $b$ grows.

\begin{proof}
We will prove the second estimate, as the first one can be obtained taking $n$ intervals
of length $\ell_i= 1$. Each element in $\Z_b[D]$ can be written as
\[
  a = \sum_{i=1}^k\left(\sum_{j=0}^{\ell_i-1}a_{ij}b^j\right)\cdot{}b^{c_i}
    = \sum_{i=1}^k A_ib^{c_i},
\]
with $0\leq A_i \leq b^{\ell_i}-1$, for each $1\leq i \leq k$. Consequently the difference
of two elements $a,a' \in \Z_b[D]$ can be written as
\[
a-a' = \sum_{i=1}^{k} (A_i-A_i')b^{c_i},
\]
where $-b^{\ell_i}+1 \leq A_i-A_i' \leq b^{\ell_i}-1$, for each $1\leq i\leq k$. In
particular the number $\Delta$ of positive differences of two elements of $\Z_b[D]$ is at most
\[
 \frac{ \prod_{i=1}^k (2b^{\ell_i}-1)-1 }{2}.
\]

For any real number $r$, let $\T_r=\R/r\Z$, and let $\pi_r : \R \rightarrow \T_r$ be the projection
map. For $z\in\Z\setminus\{0\}$ the map $\beta : \T_r \rightarrow \T_r$ given by multiplication by
$z$ is measure-preserving, i.e. for each measurable $X\subseteq \T_r$ the measure of $\beta^{-1}(X)$
is equal to the measure of $X$.

Let now $b^m$ be a power of $b$ which is $> 2\Delta$. The measure in $\T_{b^m}$ of the set
$U_z=\beta^{-1}\left(\pi_{b^m}{([-1,1])}\right)$ of the `bad' $\lambda \in \T_{b^m}$ such that
$z\lambda\in\pi_{b^m}{([-1,1])}$ is equal to $2/{b^m}$, supposing the measure of $\T_{b^m}$ to be
normalized to $1$. Since the number of positive and non-zero differences is $< {b^m}/2$, and clearly
$U_{-z} = U_z$, we have that the union of all the $U_{a-a'}$ for all distinct $a,a' \in \Z_p[D]$
cannot be all $ \T_{b^m}$.  Consequently since $\Z_{(b)}$ is dense in $\R$ there exist an element in
$\nu \in \Z_{(b)}$ which falls out of all the $U_{a-a'}$ when reduced modulo ${b^m}$, being the
above union a closed set. We have that $a\nu$ and $a'\nu$ differ by at least $1$ after reduction
modulo ${b^m}$, and hence the map $a \mapsto \dig{a\cdot \nu}{0,m}$ is injective.

Multiplying by the smallest power of $b^k$ divisible by the denominator
of $\nu$ obtain an integer $\mu = \nu b^k$ with the
required properties with respect to the map $a \mapsto \dig{a\cdot
  \mu}{k,m}$.
\end{proof}

\bibliographystyle{plain} \bibliography{biblio}

\end{document}